\documentclass[12pt]{article}
\usepackage{amsmath,amsfonts}
\usepackage[mathscr]{eucal}
\usepackage{amssymb,amsmath}
\usepackage{latexsym}
\usepackage{amsfonts}
\usepackage{amscd}
\usepackage{amsthm}
\usepackage{graphics}
\usepackage{epsfig}
\newcommand {\h}{h}
\newcommand {\side}{\textrm{side}}
\newcommand {\sside}{\textrm{sd}}
\newcommand {\del}{\partial}
\newcommand {\bt}{\tilde{\beta}}
\newcommand {\btk}{\tilde{\beta}^{(k)}}

\newcommand {\R}{\mathbb {R}}

\newcommand {\cH}{\mathcal {H}}
\newcommand {\M}{\mathcal {M}}

\newcommand {\E}{\mathcal E}
\newcommand {\U}{\mathcal U}
\newcommand {\cD}{\mathcal {D}}

\newcommand {\Z}{\mathbb {Z}}

\newcommand{\diam}{{\rm diam}}
\newcommand{\sdiam}[1]{\lvert{#1}\rvert}

\newcommand{\dist}{{\rm dist}}
\newcommand{\ball}{{\rm Ball}}

\newtheorem{theorem}{Theorem}[section]
\newtheorem{lemma}[theorem]{Lemma}

\newtheorem{cor}[theorem]{Corollary}
\newtheorem{rem}[theorem]{Remark}
\numberwithin{equation}{section}

\begin{document}

\title{Bi-Lipschitz Decomposition of Lipschitz functions into a Metric space
}

\author{Raanan Schul}

\date{}

\maketitle

\psdraft

\begin{abstract}
We prove a quantitative version of the following statement. Given a  Lipschitz function $f$ from the k-dimensional unit cube into a general metric space, one can be decomposed $f$ into a finite number of BiLipschitz functions $f|_{F_i}$ so that the k-Hausdorff content of 
\phantom{} $f([0,1]^k\smallsetminus \cup F_i)$ 
is small.  We thus generalize a theorem of P. Jones \cite{Jones-lip-bilip} from the setting of $\R^d$ to the setting of a general metric space.  This positively answers problem 11.13 in ``Fractured Fractals and Broken Dreams" by G. David and S. Semmes, 
or equivalently, question 9 from ``Thirty-three yes or no questions about mappings, measures, and metrics" by J. Heinonen and S. Semmes. 
Our statements extend to the case of {\it coarse} Lipschitz functions. 
\end{abstract}

\section{Introduction}
We prove the following theorem.
\begin{theorem}\label{t:thm-1}
Let $\epsilon\geq0$, $0<\alpha<1$ and $k\geq 1$ be  given.
There are universal constants  $M=M(\alpha,k)$, $c_1=c_1(k)$ and $c_2$ such that the following statements hold.
Let $\M$ be any metric space.
Let $f:[0,1]^k\to \M$ be an $\epsilon$-coarse 1-Lipschitz function, i.e. such that 
\begin{equation*}
\dist(f(x),f(y))\leq |x-y| +\epsilon\,.
\end{equation*}
Then there are sets 
$F_1,...,F_M\subset [0,1]^k$ so that for $1\leq i\leq M$,  
$x,y\in F_i$ we have
\begin{equation*}
\alpha|x-y| -c_2\epsilon \leq \dist(f(x),f(y))\leq |x-y|+ \epsilon\,,
\end{equation*}
and 
\begin{gather}\label{hk-estimate}
\h^k(f([0,1]^k\setminus(F_1\cup...\cup F_M)))\leq c_1 \alpha\,.
\end{gather}
($\h^k$ is the one-dimensional Hausdorff content, defined below).
\end{theorem}
As a corollary of the theorem (with $\epsilon=0$) we get 
\begin{cor}\label{HS-cor}
Let $k>0$.  Let $\M$ be a $k$-Ahlfors-David regular metric space which has Big Pieces of Lipschitz images of $\R^k$.  
Then $\M$ has  Big Pieces of BiLipshcitz images of $\R^k$.
\end{cor}
Thus we positively answer problem 11.13 in \cite{broken-dreams} by G. David and S. Semmes, or equivalently, question 9 in \cite{HS-33} by J. Heinonen and S. Semmes.

The statement of Corollary \ref{HS-cor} for the case where  $\M=\R^d$ was 
proved by Guy David  \cite{David-lip-bilip}.
In fact, David assumes less about the domain of the functions in question.  In particular, the domains need not be Euclidean, however they are required to satisfy some geometric conditions. 
Shortly after, Peter Jones gave a proof for  
Theorem \ref{t:thm-1} for the case where  $\{\epsilon=0$ and $\M=\R^d\}$  \cite{Jones-lip-bilip}.   
(Jones' and David's results appeared in the same issue of Rev. Mat. Iberoamericana.)
Their work was motivated by the study of singular integrals.

Theorem \ref {t:thm-1}, with $\epsilon=0$, can be thought of as a  quantitative version of Sard's Theorem, where we think of the 
 non-quantitative version as: a Lipschitz map can be written as a countable union of invertible maps, and a map whose range is a Lebesgue null set. 
We note that a non-quantitative variant of our theorem had already appeared in \cite{Kirchheim}.

The  proof we give  follows the outline of \cite{Jones-lip-bilip}. 
An important point is that the $\R^d$  result relies on a  sum of squares of wavelet coefficients (see  the exposition in \cite{Da})  or their upper half space analogue (which was the  way the proof went in \cite{Jones-lip-bilip}). 
We replace this by a statement about a metric space analogue of the Jones $\beta$ numbers, or equivalently, a statement about certain Menger curvature averages (Lemma \ref{sum-beta-k-leq-1}).
It is the authors feeling that the lack of such a statement was the only thing that prevented this theorem from appearing 10-20 years ago.
This theorem is another building block in the process of transferring (parts of) the Euclidean theory of quantitative rectifiability to the setting of general metric spaces.

Let us quickly define the relevant notions.
A metric space $\M$ is said to be a $k$-Ahlfors-David regular 
(with constant $c_1$) if 
for any $x\in\M,\ \ 0<r<\diam(\M)$ we have 
$c_1^{-1}r^k\leq \cH^k(\ball(x,r))\leq c_1 r^k$.
A  $k$-Ahlfors-David regular metric space $\M$ is said to  have Big Pieces of Lipschitz images 
(with constants $L_1$ and $c_1$) if 
for any $x\in\M,\ \ 0<r<\diam(\M)$ we have an $L_1$ Lipschitz function 
$f:A \to \M$ , where $A\subset \ball_{\R^k}(0,r)$
such that 
$\cH^k(f(A)\cap\ball(x,r)) \geq c_1r^k$.
A  $k$-Ahlfors-David regular metric space $\M$ is said to  have Big Pieces of BiLipschitz images 
(with constants $L_2$ and $c_2$) if 
for any $x\in\M,\ \ 0<r<\diam(\M)$ we have an $L_2$ BiLipschitz function 
$f:A \to \M$ , where $A\subset \ball_{\R^k}(0,r)$
such that 
$\cH^k(f(A)\cap\ball(x,r)) \geq c_2r^k$.

See  \cite{DS} or \cite{broken-dreams} for more details.  
We note that in question 9 of  \cite{HS-33}  there is an error in the definition of Big Pieces of Lipschitz images.
\begin{proof}[Proof of Corollary \ref{HS-cor}]
Let $c_1,\  x,\ r,\ L_1\ f,\ A$ as in the definition of Big Pieces of Lipschitz images be given.   
Let $e:\M \to L^\infty(\M)$ be the Kuratowski embedding.  Using the McShane-Whitney extension
lemma for each coordinate, we extend $e\circ f$ to  a $L_1$-Lipschitz function $\tilde f:[-r,r]^k\to L^\infty(M)$.  See page 10 in \cite{Heinonen_embedding_lec} for more details.
We now apply  Theorem \ref{t:thm-1} (rescaled) to $\tilde f:[-r,r]^k\to L^\infty(M)$ with sufficiently small $\alpha$ (depending  on the $k$-Ahlfors-David-regularity constant of $\M$, as well as $L_1$, and $k$) and $\epsilon=0$ to get 
$\cH^k (f(E\setminus(\cup F_i)))\leq \frac12 c_1 r^k $. Hence one of the sets $E\cap F_i$ must satisfy 
$\cH^k(f(E\cap F_i)\cap \ball(x,r))\geq \frac{c_1}{2M} r^k$, as desired. 
\end{proof}

\subsection*{Acknowledgments}
The author would like to thank Stephen Keith for pointing out a mistake in an earlier version of this essay. We also thank the anonymous referee for pointing out a significant oversight in the  discussion of history which follows Corollary \ref{HS-cor}.
The author is  partially  supported by NSF DMS 0502747.
\section{Proof of Theorem \ref{t:thm-1}}
\subsection{Definitions}
For a set $E$, define the one-dimensional Hausdorff content of a set $K$ as
\begin{equation*}
h^k(E)=\inf \{\sum\diam(U_i)^k: \cup U_i \supset E\}\,.
\end{equation*}
Let $p$ be a function with range contained in $\M$.
Define $\del_1(x,y,z)=\del_1^{(p)}(x,y,z)$ by
\begin{equation*}
\del_1(x,y,z)=\dist(p(x),p(y)) + \dist(p(y),p(z)) - \dist(p(x),p(z))\,.
\end{equation*}
Define for an interval $I=[a,b]\subset\R$ the quantity $\bt(I)=\bt_{(p)}(I)$ by
\begin{equation*}
\bt(I)^2\diam(I)=\diam(I)^{-3}\int_{x=a}^{x=b}\int_{y=x}^{y=b}\int_{z=y}^{z=b} \del_1(x,y,z)dzdydx\,.
\end{equation*}
We extend this definition to higher dimensional cubes  by rotations.  Let $k>1$.  
Define for a cube $Q\in \R^k$ the quantity $\btk(Q)=\btk_{(p)}(Q)$ by
\begin{eqnarray*}
&&\btk(Q)^2\side(Q)^{k-1}=\\
&&\phantom{xxx}
\int_{g\in G_k}\int_{x\in \R^k \circleddash g\R} \chi_{\{|(x+g\R)\cap 7Q|\geq \side(Q)\}}\bt((x+g\R)\cap 7Q)^2dxd\mu(g)
\end{eqnarray*}
where  
$\R$ is identified with $\{\R,0,...,0\}\subset \R^k$,  $G_k$ is the group of all rotations of $\R$ in $\R^k$ equipped with the its Haar measure $d\mu$, and $dx$ is the $k-1$ dimensional Lebesgue measure on $\R^k \circleddash g\R$, the orthogonal complement of $g\R$ in $\R^k$.
We write $\bt^{(1)}=\bt$.
and note that  any $k\geq 1$, we have that  $\btk$ is scale invariant.
This type of quantity is connected to Menger curvature.  See \cite{RS-TSP-survey} for more details.

\noindent
Define $\cD_0$ the standard dyadic partition of $\R$, i.e.
\begin{equation*}
\cD_0:=\{ [\frac{j_1}{2^{j_2}}, \frac{j_1+1}{2^{j_2}}]:j_1,j_2\in\Z\}\,.
\end{equation*}
Define $\cD_1$ a dyadic partition of $\R$ given by shifting the standard dyadic partition by $\frac13$, i.e. 
$\cD_0 + \frac13$.
For $i=(i_1,...,i_k)\in\{0,1\}^k$ we define
\begin{equation*}
\cD^k_i:=\cD_{i_1}\oplus....\oplus\cD_{i_k}\,.
\end{equation*}
The fact that now a ball $\ball(x,r)\subset \R^k$ with $r< \frac16$ is contained in a cube $Q\in \bigcup\limits_{i\in\{0,1\}^k} \cD^k_i $ with $\side(Q)\sim r$ earns this setup the (now standard) name {\it the one third trick}. 

For simplicity of notation in the proof, we extend $f$ to be 1-Lipschitz with domain  $\R^k$ (by say fixing $f$ on rays emanating from $(\frac12,...,\frac12)$ and outside $(0,1)^k$. 

\noindent
We call two dyadic cubes  $Q_1$ and $Q_2$ {\bf semi-adjacent} if 
\begin{equation*}
0 < \dist(Q_1,Q_2)\leq 2\diam(Q_1)=2\diam(Q_2)\,.
\end{equation*}
Hence every cube $Q$ has at most $C(k)$ semi-adjacent cubes.

\subsection{The proof}
We start by defining a Lipschitz function $p$.

Let $X=X_{\epsilon}\subset [0,1]^k$ be an $\epsilon$-net for $[0,1]^k$ and 
$Z=Z_{\epsilon}\subset \M$ be an $\epsilon$-net for $\M$.
Consider a function $f':X\to \M$ defined as follows.
for $z\in \M$, let $z_\epsilon\in Z$ be such that $\dist(z,z_\epsilon)\leq \epsilon$ (chosen arbitrarily if there is more then one such $z_\epsilon$). 
Define $f'(x)=f(x)_\epsilon$. 
We have for any $x,y\in X$ 
\begin{equation*}
\dist(f'(x),f(x))\leq \epsilon,\quad
\dist(f'(x),f'(y))\leq |x-y| + \epsilon +2\epsilon\,,
\end{equation*}
We get that  for any $x,y\in X$
\begin{equation*}
\dist(f'(x),f'(y))\leq 4|x-y|\,.
\end{equation*}
Now extend $f'$ to a 4-Lipschitz function $p:[0,1]^k\to L^\infty(\M)$ using 
the Kuratowski embedding and  the McShane-Whitney extension
(as in Corollary  \ref{HS-cor}).
Denote  by $\tilde f:[0,1]^k\to L^\infty(\M)$,   the map given by  using the Kuratowski embedding of $\M$ in $L^\infty(\M)$.
We have for $x\in [0,1]^k$ and $x_\epsilon\in X$ such that $\dist(x,x_\epsilon)\leq \epsilon$
\begin{equation}\label{p-and-f-close}
\begin{aligned}
\dist( \tilde f(x),p(x))&\leq
\dist( \tilde f(x),\tilde f(x_\epsilon)) +
\dist( \tilde f(x_\epsilon),p(x_\epsilon))+
\dist( p(x_\epsilon),p(x))\\
&\leq
2\epsilon + \epsilon + 4\epsilon\,.
\end{aligned}
\end{equation}

This $p$ is the one we use in the above definitions of  $\del_1$ and  $\beta$.

\begin{lemma}\label{sum-beta-k-leq-1}
For an L-Lipschitz function $p$,
\begin{equation*}
\sum_{Q\in \cD^k_i,\ Q\subset [0,1]^k \atop i\in\{0,1\}^k} \btk(Q)^2\side(Q)^k  \lesssim L\,.
\end{equation*}
\end{lemma}
We postpone the proof of this lemma to Section \ref{curvature_estimates}.

Let $\alpha'=10\alpha$.
For $x_1,x_2\in \R^k$, let $[x_1,x_2]$ be the straight segment connecting $x_1$ and $x_2$. 
Let 
\begin{eqnarray*}
\E_1&:=\Big\{Q_1\in\cD^k_{\bar{0}}:&Q_1\subset [0,1]^k,\quad \exists x_1\in Q_1,\ x_2\in Q_2, \\  
	&& Q_1,Q_2\textrm{ semi-adjacent}, \ \  
	       \diam(p([x_1,x_2]))\geq {\alpha'} |x_1-x_2|  , \\
	&&  \dist(p(x_1),p(x_2))\leq \frac{\alpha'}{10} |x_1-x_2|, \quad {\alpha'}|x_1-x_2|\geq 10\epsilon  \Big\}\,,\\
\E_2&:=\Big\{[x_1,x_2]:&x_i\in Q_i, \quad Q_1\subset [0,1]^k,\quad\\
			&&Q_1,Q_2\textrm{ semi-adjacent}, \quad 
			\diam(p([x_1,x_2]))\leq {\alpha'} |x_1-x_2|\,, \\
			&&   {\alpha'}|x_1-x_2|\geq 10\epsilon\Big\}\,.\\
B&:=\Big\{Q_0\in \E_1:&\exists Q_1,...,Q_N\in\E_1, \textrm{ such that } 
	Q_0\subsetneq Q_1\subsetneq...\subsetneq Q_N\Big\}\,.
\end{eqnarray*}
The constant $N$ will be chosen later, and will depend only on ${\alpha'}$ and $k$.
(Note that in the definition of $B$ the cubes $Q_1,...,Q_N$ may be of wildly different scales.)

\begin{lemma}
\begin{equation*}
\h^k(f( \cup\E_2))\lesssim {\alpha'}\,.
\end{equation*}
\end{lemma}	
\begin{proof}
Assume $[x_1,x_2]\in \E_2$, and let $Q_1, Q_2$ be the corresponding semi-adjacent cubes.
Recall that $|x_1-x_2|\sim \diam(Q_1)$.
Define for $c=3,30$
\begin{equation*}
U^c_{x_1,x_2}=f^{-1}\ball(f(x_1),c{\alpha'}|x_1-x_2|)\,.
\end{equation*}
Then 
\begin{equation}\label{large-u2}
U^3_{x_1,x_2}\supset \{x\in \R^k:\dist(x,[x_1,x_2])<{\alpha'} |x_1-x_2|\}\,.
\end{equation}
Note that
\begin{equation*}
\diam( fU^{30}_{x_1,x_2} )\lesssim {\alpha'} |x_2-x_1| 
\end{equation*}
implying that (together with \eqref{large-u2})
\begin{equation}\label{small-diam}
\diam( fU^{30}_{x_1,x_2} )^k\lesssim {\alpha'} \cH^k(U^3_{x_1,x_2})\,.
\end{equation}
Consider the set 
\begin{equation*}
\U=\{U^3_{x_1,x_2}: [x_1,x_2]\in \E_2\}\,.
\end{equation*}
We will show $h^k(\cup f\U)\lesssim {\alpha'}$, which will give the lemma as $\cup \E_2\subset \cup\U$.

We use a Vitali covering type argument.
We  find a disjoint sub-collection $\U'\subset \U$ so that if 
$U\in \U$ then $U\cap U'\neq \emptyset$ for some  $U'\in \U'$ with $2\diam(fU')\geq\diam(fU)$ as follows. 
Write $\U=\cup \U_j$ where $U\in \U_j$ implies
\begin{equation*}
2^{-j-1}<\diam(fU)\leq 2^{-j}\,.
\end{equation*}
We greedily construct $\U'$ by adding sets to it from $\U_j$, inducting on $j$.
We start with $\U'=\emptyset$.
Place a maximal (with respect to inclusion) disjoint subset of $\U_0$ in $\U'$.
At stage $j>0$, consider all sets $S_j\subset \U_j$ which have  disjoint elements and have elements disjoint from all current $\U'$ elements. Take a maximal (with respect to inclusion) such $S_j$, and add $S_j$ to $\U'$. This defines $\U'$ as desired.

Now, we note that if $U\in \U$ and $U\cap U'\neq \emptyset$, $U'\in \U'$ 
with $2\diam(fU')\geq\diam(fU)$, and $U'=U^3_{x_1,x_2}$, then (by looking at the push-forward by $f$) 
$U\subset U^{30}_{x_1,x_2}$.
Hence 
\begin{equation*}
\cup\{U:U\in \U\} \subset \cup\{U^{30}_{x_1,x_2}:U^3_{x_1,x_2}\in \U'\}\,.
\end{equation*}
Using  the disjointness of elements in $\U'$  and 
inequality \eqref{small-diam}
we get the desired result.
\end{proof}

\begin{lemma}\label{large-beta}
There is an $\epsilon_0=\epsilon_0({\alpha'},k)>0$ such that 
for any  $Q_1\in\E_1$  we have 
 $\btk(7Q_1)\geq \epsilon_0$.  
\end{lemma}
\begin{proof}
Let $Q_2,x_1,x_2$ be as in the definition of $\E_1$. Set $D=\side(Q_1)\geq |x_2-x_1|$.
Let $x_3\in [x_1,x_2]$ be a point such that $\dist(p(x_3),p(x_1))> \frac{{\alpha'}}2D$.
Then   for $x_i'\in\ball(x_i,\frac{\alpha'}{100} D)$ we have 
$\del_1(x_1',x_3',x_2')\geq \frac{{\alpha'}}{10}D$.
Now use definition of $\btk$.
\end{proof}

Let $B_1=\cup_{Q\in B\atop Q\subset [0,1]^k} 7Q$.
From Lemma \ref{sum-beta-k-leq-1}, Lemma \ref{large-beta}, and the one-third trick,  we have:
\begin{lemma}
By taking $N=N({\alpha'},k)$ large enough (universally determined)
\begin{equation*}
\cH^k(B_1)\lesssim {\alpha'}\,.
\end{equation*}
and
\begin{equation*}
\h^k (f(B_1))\lesssim {\alpha'}\,.
\end{equation*}
\end{lemma}
\begin{proof}
The first inequality follows  from 
\begin{equation*}
\int \chi_{B_1}
\leq 
		\int \frac1{N+1}\sum\limits_{Q\in \E_1} \chi_{7Q}\leq 
		\frac{C}{N+1} \,.
\end{equation*}
If $\epsilon=0$ this is more then enough for the second inequality as well.
If $\epsilon>0$, to see the second inequality we note that
\begin{eqnarray*}
\h^k (f([0,1]^k\cap B_1))&\leq&
\frac1{N+1}\sum\limits_{Q\in \E_1} \side(7Q)+\epsilon\\
&\leq& 
\frac1{N+1}\sum\limits_{Q\in \E_1} 2\side(7Q) \leq 
		\frac{C}{N+1} \,.
\end{eqnarray*} 
\end{proof}

Now denote by $G=[0,1]^k\setminus (B_1 \cup (\cup\E_2))$. We would like to split $G$ into $M({\alpha'},k)$ sets as desired.  
We split according to the behavior of the function $p$ using $\E_1$ as our guide.
One goes through the dyadic tree (large scale  to fine scale) and makes sure that if $Q_1,Q_2\in\E_1$ are semi-adjacent, then they are not in the same $F_i$.  
Since we excise the intervals in $B$ (which gave us $B_1$), 
we can do this with only a finite number of sets $F_i$ 
(namley, $M=2^{C(k)N}$, with $N\sim\frac1{\alpha'}$).
For more details see
pages 81-82  of \cite{Da} (starting at the bottom of page 81, with the same notation). 
This gives
$F_1,...,F_M\subset [0,1]^k$ so that for $1\leq i\leq M$,  
and $x,y\in F_i$ such that $\alpha'|x_1-x_2|\geq 10\epsilon$ we have
\begin{equation*}
\frac1{10}{\alpha'}|x-y| \leq \dist(p(x),p(y))\,,
\end{equation*}
and 
\begin{gather*}
\h^k(f([0,1]^k\setminus(F_1\cup...\cup F_M)))\leq c_1 {\alpha'}\,.
\end{gather*}
Using equation \eqref{p-and-f-close} we get that if
$x,y\in F_i$ we have
\begin{equation*}
\frac1{10}{\alpha'}|x-y| -14\epsilon \leq \dist(f(x),f(y))\,,
\end{equation*}
This concludes the proof of Theorem \ref{t:thm-1}.

\section{Curvature estimates}\label{curvature_estimates} 
In this section we prove  Lemma \ref{sum-beta-k-leq-1}.
We first consider the case $k=1$, and then use it to prove the lemma for $k>1$.
\begin{lemma}\label{sum-beta-leq-1}
\begin{equation*}
\sum_{\cD_0\cup\cD_1\atop  I\subset [0,1]} \bt(I)^2\diam(I)  \lesssim L\,.
\end{equation*}
\end{lemma}
This lemma is stated and proved in \cite{RS-metric}. 
The setting we were interested in there was that of {\it Ahlfors-regular curves}, however the proof given there for this lemma is correct for the setting we have here.  
We give most of the proof's details 
in the appendix of this paper.

We are now ready to prove Lemma \ref{sum-beta-k-leq-1}.
\begin{proof}[Proof of Lemma \ref{sum-beta-k-leq-1}]
We use Lemma \ref{sum-beta-leq-1} and the definition of $\btk$.
Fix $i\in\{0,1\}^k$  and write $\cD^k=\cD^k_i$. 
We have
\begin{eqnarray*}
&&\sum_{Q\in \cD^k \atop Q\subset [0,1]^k} \btk(Q)^2\sside(Q)^k\\
&=&
\sum_{Q\in \cD^k \atop Q\subset [0,1]^k} 
	\int\limits_{\ g\in G_k}\int\limits_{\ x\in \R^k \circleddash g\R} \chi_{\{|(x+g\R)\cap 7Q|\geq \sside(Q)\}}\bt((x+g\R)\cap 7Q)^2\sside(Q)
							dxd\mu(g)\\
&=&
\int\limits_{\ g\in G_k}\int\limits_{\ x\in \R^k \circleddash g\R} 
	\sum_{Q\in \cD^k \atop Q\subset [0,1]^k} \chi_{\{|(x+g\R)\cap 7Q|\geq \sside(Q)\}}\bt((x+g\R)\cap 7Q)^2\sside(Q)
							dxd\mu(g)\\
&\lesssim&
L+\int\limits_{\ g\in G_k}\int\limits_{\ x\in \R^k \circleddash g\R} 
	\chi_{\{x\in C[0,1]^k\}} \sum_{I\in \cD_0\cup\cD_1\atop I\subset [0,1]} \bt(x+gI)^2\diam(I)
							dxd\mu(g)\\
&\lesssim&
L+ L\int\limits_{\ g\in G_k}\int_{x\in \R^k \circleddash g\R} 
	\chi_{\{x\in C[0,1]^k\}} 	dxd\mu(g) 
	\lesssim
L+ L\int\limits_{\ g\in G_k}d\mu(g)
\lesssim
L\,.
\end{eqnarray*}
where above, the notation $\sside(Q)$ is short for $\side(Q)$, the side length of the cube $Q$. 
\end{proof}

\section{Appendix}
We review the proof of Lemma \ref{sum-beta-leq-1}, taken from \cite{RS-metric}. For a little more details see the original.

For  numbers $r,v\in[0,1]$  we will look at the mapping 
$\psi^{v,r}:[0,1]\to[0,1]$ given by $\psi^{v,r}(t)= v+ r t \mod 1$. 

For an interval $I\subset[0,1]$ write $I=[a(I),b(I)]$. 
\begin{rem}
Let $\sdiam{I}$ be the diameter of the interval $I$.
When doing addition  $\mod 1$, we have (by change of variable) for any $I'$ with $\sdiam{I'}=2^{-k}$ 
\begin{eqnarray*}
&&\sum_{I\in\cD_0\atop \sdiam{I}=2^{-k}}
\sdiam{I}^{-3}
	\int\limits_{a(I)}^{b(I)}\ \int\limits_x^{b(I)}\ \int\limits_y^{b(I)}
		\del_1(x,y,z)dzdydx \\
&\leq&
\sdiam{I'}^{-3}\int\limits_{v=0}^1\ 
	\int\limits_{r=0}^1\ 
			\int\limits_{y\in v+rI'}
				\del_1(v+ra(I'),y,v+rb(I')) 
					 dy \cdot \sdiam{I'}dr dv
\end{eqnarray*}
giving 
\begin{eqnarray*}
&&\sum_{I\in\cD_0\atop \sdiam{I}=2^{-k}}
\sdiam{I}^{-3}
	\int\limits_{a(I)}^{b(I)}\ \int\limits_x^{b(I)}\ \int\limits_y^{b(I)}
		\del_1(x,y,z)dzdydx \\
&\leq&
\sum_{I\in\cD_0\atop \sdiam{I}=2^{-k}}
\sdiam{I}^{-3}\int\limits_{v=0}^1\ 
	\int\limits_{r=0}^1\ 
			\int\limits_{y\in v+rI}
				\del_1(v+ra(I),y,v+rb(I)) 
					 dy \cdot \sdiam{I}dr  \sdiam{I}dv\,.
\end{eqnarray*}
\end{rem}

Let $I'=[a,b]\in\cD_0$.
Define
\begin{gather*}
\del_{dyadic}( \psi^{v,r}(I')) := \del_1(		(v+ra),
									(v+r{a+b\over 2}),
									(v+rb)	).
\end{gather*} 

The triangle inequality gives the following lemma.
\begin{lemma}
Let $I\in\cD_0$.
Let $v,r\in[0,1]$ be chosen such that $\psi^{v,r}(I)=[x,z]\ni y$. Then
\begin{gather*}
\del_1(x,y,z) \leq 
\sum\limits_{I'\in\cD, I'\subset I \atop y\in \psi^{v,r}(I')}
	\del_{dyadic}(\psi^{v,r}(I')).
\end{gather*}
\end{lemma}

Via telescoping sums, one get the following lemma.
\begin{lemma}
Let $r,v\in[0,1]$ be fixed.
Then
\begin{gather*}
\sum\limits_{I'\in\cD_0}
	\del_{dyadic}(\psi^{v,r}(I')) 
\lesssim
L
\end{gather*}
\end{lemma}
Finally, putting the above together, one gets (with the diameter of $I$ denoted by $\sdiam{I}$)
\begin{eqnarray*}
&&\sum_{I\in\cD_0}
\sdiam{I}^{-3}
	\int\limits_{a(I)}^{b(I)}\ \int\limits_x^{b(I)}\ \int\limits_y^{b(I)}
		\del_1(x,y,z)dzdydx \\
&\leq&
\sum_{I\in\cD_0}
\sdiam{I}^{-3}\int\limits_{v=0}^1\ 
	\int\limits_{r=0}^1\ 
			\int\limits_{y\in v+rI}
				\del_1(v+ra(I),y,v+rb(I)) 
					 dy \cdot \sdiam{I}dr \cdot \sdiam{I}dv\\
&\leq&
\sum_{I\in\cD_0}
\sdiam{I}^{-3}\int\limits_{v=0}^1\ 
	\int\limits_{r=0}^1\ 
		\sum\limits_{I'\in\cD_0\atop I'\subset I}\ 
			\int\limits_{y\in v+rI'}
				\del_{dyadic}(\psi^{v,r}(I')) 
					\cdot  dy \cdot \sdiam{I}dr \cdot \sdiam{I}dv\\
&=&
\sum_{I\in\cD_0}
\sdiam{I}^{-3}\int\limits_{v=0}^1\ 
	\int\limits_{r=0}^1\ 
		\sum\limits_{I'\in\cD_0\atop I'\subset I}
				\del_{dyadic}(\psi^{v,r}(I')) 
					\cdot r \cH^1(I') \cdot \sdiam{I}dr \cdot \sdiam{I}dv\\
&=&
\int\limits_{v=0}^1\ 
	\int\limits_{r=0}^1\ 
	\sum_{I\in\cD_0}
		{1\over \sdiam{I}} 
			\sum\limits_{I'\in\cD_0\atop I'\subset I}
				\del_{dyadic}(\psi^{v,r}(I')) 
					\cdot r\cH^1(I')  dr dv\\
&=&
\int\limits_{v=0}^1\ 
	\int\limits_{r=0}^1\ 
		\sum_{I'\in\cD_0}\ 
				\sum\limits_{I\supset I'}
					{\cH^1(I')\over \sdiam{I}} 
						\del_{dyadic}(\psi^{v,r}(I'))
						\cdot r  dr dv\\
&\lesssim&
\int\limits_{v=0}^1\ 
	\int\limits_{r=0}^1\ 
		\sum_{I'\in\cD_0} 
						\del_{dyadic}(\psi^{v,r}(I'))
						\cdot r  dr dv
\lesssim L.
\end{eqnarray*}
which gives Lemma \ref{sum-beta-leq-1}.

\bibliographystyle{alpha}
\bibliography{../../bibliography/bib-file-1}   

\bigskip
\section*{}
{ \flushright
Raanan Schul\\
UCLA Mathematics Department\\
Box 951555\\
Los Angeles, CA 90095-1555\\
U.S.A.\\
{\tt schul@math.ucla.edu}\\
}

\end{document}